\documentclass{amsart}
\usepackage{amssymb,stmaryrd}
\usepackage{amsfonts}
\usepackage{amstext}
\usepackage{algorithmic}
\usepackage{algorithm}
\usepackage{graphicx}
\usepackage{epstopdf}
\usepackage[all]{xy}
\usepackage{MnSymbol}

\numberwithin{equation}{section}
\newtheorem{theorem}{Theorem}[section]

\newtheorem{proposition}[theorem]{Proposition}

\theoremstyle{definition}

\theoremstyle{remark}

\newcommand{\R}{{\mathbb{R}}}

\newcommand{\C}{{\mathbb{C}}}

\newcommand{\pdol}{\overline{\partial}}
\newcommand{\pd}{\partial}

\newcommand{\CR}{{\mathcal{R}}}

\newcommand{\wedgeq}{{\wedge\kern-5pt\cdot}}

\newcommand{\tens}{\otimes}

\newcommand{\id}{{\rm id}}

\newcommand{\extd}{{\rm d}}
\newcommand{\del}{{\partial}}

\begin{document}

\title{Quantum Riemannian geometry of phase space and nonassociativity} 

\keywords{noncommutative geometry, quantum gravity, Poisson geometry, Riemannian geometry, quantum mechanics}


\author{ Edwin J. Beggs \& Shahn Majid}
\address{Department of Mathematics, Swansea University \\ Singleton Parc, Swansea SA2 8PP\\
+ \\
Queen Mary, University of London\\
School of Mathematics, Mile End Rd, London E1 4NS, UK}

\email{s.majid@qmul.ac.uk}
\thanks{The 2nd author was on leave at the Mathematical Institute, Oxford}

\begin{abstract}  Noncommutative or `quantum' differential geometry has emerged in recent years as a process for quantizing not only a classical space into  a noncommutative algebra (as familiar in quantum mechanics) but also differential forms, bundles and Riemannian structures at this level. The data for the algebra quantisation is a classical Poisson bracket, the data for the quantum differential forms is a Poisson-compatible connection it was recently shown in \cite{BegMa5} that after this, classical data such as classical bundles, metrics etc.\  all become quantised in a canonical `functorial' way at least to 1st order in deformation theory. There are, however, fresh compatibility conditions between the classical Riemannian and the Poisson structures as well as new physics such as nonassociativity at 2nd order. We give an introduction to this theory from \cite{BegMa5} and some details for the case of $\mathbb{C P}^n$ where the commutation relations have the canonical form $[w^i,\bar w^j]=\mathrm{i}\lambda\delta_{ij}$ similar to the proposal of Penrose for quantum twistor space. Our work provides a canonical but ultimately nonassociative differential calculus on this algebra and quantises the metric and Levi-Civita connection at lowest order in $\lambda$. \end{abstract}
\maketitle 

\section{Introduction}

There are today lots of reasons to think that spacetime itself is better modelled as  `quantum' due to Planck-scale corrections. Here the quantisation parameter $\lambda$ is the Planck scale around $10^{-33}$cm (more precisely, $\mathrm{i}$ times this as we work with imaginary $\lambda$) and in this context we should quantise not only the spacetime variables but classical Riemannian geometry or `gravity' variables as well. This is a viable approach to the hot topic of quantum gravity phenomenology as it allows one to come rather quickly to predictions for Planck scale effects. In other contexts one might have $\lambda=\imath\hbar$ and be  interested in Quantum Mechanics on a classical phase space that has other classical geometrical data on it, including perhaps a Riemannian structure, and one could ask how does this structure extend to the quantum algebra. This would be relevant for example to quantisation of the geometrical description\cite{Brody} of Berry phase. These and many other contexts where one may want for mathematical or physical reasons to `follow' classical geometry into the quantum domain can now be addressed using noncommutative geometry.  

The approach to noncommutative geometry that we use is one that has developed over the years mainly from experience with quantum groups, see in particular our papers \cite{BegMa1}-\cite{BegMa5} and references therein.  This approach is very  different from  and has a completely opposite starting point to the approach of Alain Connes\cite{Con}. The latter starts `top down' with a spectral triple as an algebraic model of the Dirac operator on spinors whereas our approach is `bottom up' starting with differential structures and ultimately, we hope, building up to spinors and a Dirac operator as a final layer that is not yet fully understood. Although our approach to noncommutative differential and Riemannian geometry now exists as a noncommutative algebraic framework, and has some fully worked examples such as the 2D analysis in \cite{BegMa4}, there still remains the general problem of construction from classical data. This was recently solved in \cite{BegMa5} as follows. 

The very first layer of the problem from our point of view is of course the Poisson structure, a tenet of mathematical physics since the early works of Dirac  being to `quantise' this to a noncommutative algebra $A_\lambda$. Let us recall that the mathematical background to this is to consider an algebra $A_\lambda$ where $\lambda$ is a formal parameter such that $A_0$ is commutative, we denote the product of the latter by omission, and
\[ a\bullet_\lambda b= ab + O(\lambda).\]
We assume that expressions can be expanded in $\lambda$ and equated order by order. In this case 
\[ a\bullet_\lambda b-b\bullet_\lambda a=\lambda\{a,b\} + O(\lambda^2)\]
defines a map $\{\ ,\ \}$ and the assumption of an associative algebra quickly leads to the necessary feature that this is a Lie bracket and the Hamiltonian vector field $\hat a:=\{a, \}$ is a derivation on $A_0$, making $A_0$ a Poisson algebra. The converse to this is the `{\em quantisation problem}': given a smooth manifold $M$ and a Poisson bracket on it, can one deform $C^\infty(M)$ to a noncommutative algebra $A_\lambda=C^\infty(M)[[\lambda]]$ as a vector space (i.e. complexifying and working over the ring of formal power series $\C[[\lambda]]$ such that the above holds. In 1994 Fedosov\cite{FedBook} gave a geometrical solution for the symplectic case where $\{\ ,\ \}$ is nondegenerate. This uses Weyl bundles over the spacetime and a flat connection to globalise the Heisenberg-type algebra associated to the symplectic structure. In 2003 Kontsevich gave a rather different solution using a sum over graphs and bidifferential operators associated to them, for any Poisson manifold. Our first innovation in \cite{BegMa5} is instead of working over the ring $\C[[\lambda]]$ we work over the ring $\C[\lambda]/(\lambda^2)$ where we formally set $\lambda^2=0$, which we call {\em semiquantisation}. Both rings are mathematical tricks: in physical applications one wants $\lambda$ to be an actual (imaginary) number meaning on the one hand for powerseries to converge and on the other hand, in our case, for $O(\lambda^2)$ terms to be physically neglectable. This should be reasonable when $\lambda$ is the Planck scale as envisaged in many (but not the only) applications; it will be hard enough to observe these order $\lambda$ corrections and corrections beyond that are likely to be undetectable and  irrelevant to current tests of quantum gravity. At this level the semiquantisation presents no problem and does not even need $\{\ ,\ \}$ to obey the Jacobi identity. On the other hand letting go of this would entail $A_\lambda$ being nonassociative when  $\lambda^2$ is considered, which we prefer to avoid. 

The second layer of the problem is to construct not only an algebra $A_\lambda$ (specifying the algebra is roughly speaking like specifying a topological space) but a `differential graded algebra' (DGA) 
 \[ \Omega(A_\lambda)=\oplus_n\Omega^n(A_\lambda),\quad \extd:\Omega^n(A_\lambda)\to\Omega^{n+1}(A_\lambda)\]
 obeying $\extd^2=0$ and the graded-Leibniz rule. This plays the role of the algebra of differential forms and is like specifying a differential structure on a space. The data for the differential structure at the semiclassical level was analysed in \cite{Haw, BegMa1} by looking at
\[ a\bullet_\lambda \extd b-( \extd b)\bullet_\lambda a=\lambda \nabla_{\hat a}\extd b + O(\lambda^2).\]
The assumption of an associative $\Omega(A_\lambda)$ and the Leibniz rule for $\extd$ requires at order $\lambda$ that 
\[ \nabla_{\hat a}(b\extd c)=\{a,b\}\extd c+ b\nabla_{\hat a}\extd c\]
\begin{equation}\label{poicomp} \extd\{a,b\}=\nabla_{\hat a}\extd b- \nabla_{\hat b}\extd a\end{equation}
(these follow easily from $[a, b\extd c]=[a,b]\extd c+ b[a,\extd c]$ and $\extd [a,b]=[\extd a,b]+[a,\extd b]$). The first requirement says that $\nabla$ is a covariant derivative along Hamiltonian vector fields $\hat a$ and the second is a Poisson-compatibility. For simplicity we will speak of a connection $\nabla_i$ in our coordinate basis but if the Poisson tensor in these coordinates is $\omega^{ij}$ then we are only really making use of the combination $\omega^{is}\nabla_s$ in all that follows, which is to say a partial connection in the case where $\omega$ is degenerate. Simply put, the semiclassical data for the quantum differential structure is a Poisson-compatible (partial) connection $\nabla$.

This brings us to the following two quantisation problems given a manifold $M$ equipped with data $(\omega,\nabla)$ as above:

\begin{quote}{\bf Problem 1}: can we quantise the data to an associative DGA $\Omega(A_\lambda)$?
\end{quote}

\begin{quote}{\bf Problem 2}: can we similarly quantise other classical geometrical structures such as a metric and its Levi-Civita connection?
\end{quote}

Problem 1 completes the `second layer' to include higher differential forms and Problem 2 represents a `third layer' to the quantisation of the geometry.  The work  \cite{BegMa5} has answered both questions in the affirmative, but only at order $\lambda$, i.e. the semiquantisation problem is now fully solved.  Problem~1 has a canonical solution at this order without needing any new data and Problem 2 also has a canonical solution when it exists, but for existence there are new equations of constraint not seen before in physics between the Poisson bracket, the Poisson connection and the classical Riemannian structure. Let $g$ be the Riemannian metric and $S$ the contorsion tensor of $\nabla$ (so that $\nabla+S$ is the Levi-Civita connection) and we let $\CR$ be a certain `generalized Ricci 2-from' which we build by contraction of $\omega$ with the curvature $R$ and torsion of $\nabla$. Then the new conditions we find are\cite{BegMa5}
\begin{equation}\label{nabla g} g_{mn;k}=0\end{equation}
\begin{equation}\label{levicondition} \CR_{mn\hat;k}-\omega^{ij}\,g_{rs}\,S^s_{jn}(R^r{}_{mki}+S^r_{km;i})+ \omega^{ij}\,g_{rs}\,S^s_{jm}(R^r{}_{nki}+S^r_{kn;i})=0\end{equation}
where the first condition ensures centrality of the quantum metric and the second ensures quantum metric compatibility of our quantum Levi-Civita connection. Here $\hat;$ is with respect to the Levi-Civita connection and $;$ is with respect to $\nabla$. We will give more details from \cite{BegMa5} in the next section.  What is significant is for the first time to have differential constraints relating the Riemann and Poisson structures on the classical manifold, as a condition for quantisability. These conditions can be severe, for example in 2D\cite{BegMa4} for a certain well-known quantum spacetime the condition (\ref{nabla g}) forces the classical metric to either have a very strong gravitational source at the origin or to correspond to an expanding cosmology depending on the sign of a parameter. 

Finally, it is already known from \cite{BegMa1} that the quantisation even of $\Omega^1(A_\lambda)$ at the next, order $\lambda^2$, level, has a obstruction the Riemann curvature of $\nabla$. We could require $\nabla$ to have zero curvature but if we are geometers this feels like restricting ourselves for no reason to `flatland', particularly if we take the simplest case where $\nabla=\widehat\nabla$. The alternative is that we must allow nonassociativity of differential forms with functions as second order:
\[ (a\bullet_\lambda\extd b)\bullet_\lambda c- a\bullet_\lambda((\extd b)\bullet_\lambda c)=O(\lambda^2)\]
for generic functions $a,b,c$ on phase space, as controlled by the curvature. In the case where $\nabla=\hat\nabla$ we see that {\em gravity induces nonassociativity}. In this nonassociative case we see that the axioms for $\Omega(A_\lambda)$ are weaker than we had first posed. This suggests to tackle the full problem order by order in some kind of $A_\infty$ approach which remains to be worked out. The alternative to this `anomaly for differential calculus' (or Beggs-Majid no-go theorem) in \cite{BegMa1} is to absorb the anomaly by adding one or more extra dimensions to $\Omega(A_\lambda)$ which is a very different analysis\cite{Ma:rec}.

One example of classical data and where our conditions automatically hold\cite{BegMa5} is any K\"ahler-Einstein manifold. Here we take $\nabla$ to be the Levi-Civita connection so $S=0$ and $\CR$ comes out to be the usual Ricci 2-form which is covariantly constant so both (\ref{nabla g})-(\ref{levicondition}) are solved. This includes $\mathbb{C P}^n$ and we give details for semiquantisation of the differential geometry in this case. As noted by Penrose in the twistor case the quantum algebra here can be put in a canonical commutation relations form. We will describe our nonassociative calculus in these coordinates as well as in other more natural $z,\bar z$ coordinates. The `noncommutative complex structure' quantising the classical one of $\mathbb{C P}^n$ will also be touched upon, with details to appear in \cite{BegMa6}. 

\section{Semiquantum Riemannian geometry\cite{BegMa5}}

We will work in a local coordinate basis for our manifold $M$ so that 
\[ \{a,b\}=\omega^{ij}a_{,i}b_{,j},\quad \nabla_j\extd x^i=-\Gamma^i_{jk}\extd x^k\]
 where $\Gamma$ are the Christoffel symbols of our (partial) linear connection $\nabla$. Then (\ref{poicomp}) can be written \cite{BegMa5} as
\begin{equation}\label{T1} \omega^{ij}{}_{;m}+\omega^{ik}T^j_{km}-\omega^{jk}T^i_{km}=0\end{equation}
where $T^i_{jk}=\Gamma^i{}_{jk}-\Gamma^i{}_{kj}$ is the torsion tensor of $\nabla$ and $;$ is with respect to $\nabla$. Given this identity, the requirement that
the antisymmetric bivector $\omega^{ij}$ defines a Poisson tensor becomes\cite{BegMa5}
\begin{equation}\label{T2} \sum_{{\rm cyclic}(i,j,k)}\omega^{im}\omega^{jn}T^k_{mn}=0\end{equation}
again depending only on the torsion. We will speak throughout about a Poisson tensor $\omega$ and $\nabla$ Poisson-compatible, since we ultimately prefer $A_\lambda$ to be associative at all orders, but in fact we will never make use of (\ref{T2}) in what follows. 

\subsection{Quantisation of the exterior algebra}

 Given such data $(\omega,\nabla)$, the obvious structure of $\Omega^1(A_\lambda)$ at order $\lambda$ is
\[ a\bullet_\lambda b=ab+{\lambda\over 2}\omega^{ij}a_{,i}b_{,j}\]
\[ a\bullet_\lambda\xi=a\xi+{\lambda\over 2}\omega^{ij}a_{,i}\nabla_j \xi,\quad \xi\bullet_\lambda a=a\xi-{\lambda\over 2}\omega^{ij}a_{,i}\nabla_j \xi\]
for all $a,b\in C^\infty(M)$ and $\xi\in\Omega^1(M)$. This is as in \cite{BegMa1} and we extend this now to all degrees:

\begin{theorem} \cite{BegMa5} The above data extends at order $\lambda$ to a DGA $\Omega(A_\lambda)$ quantising the exterior algebra $\Omega(M)$. 
\end{theorem}
Here the quantum wedge product has a functorial part $\wedge_Q$ and a `quantum' correction
\[ \xi\wedge_1\eta=\xi\wedge_Q\eta+\lambda\,(-1)^{|\xi|+1}\, H^{ij}\wedge (\partial_i \, \righthalfcup\, \xi)\wedge(\partial_j \, \righthalfcup\, \eta);\quad  \xi\wedge_Q\eta=\xi\wedge\eta+{\lambda\over 2}\omega^{ij}\nabla_i\xi\wedge\nabla_j\eta\]
where the quantum correction is controlled by a family of 2-forms
\begin{equation}\label{H} H^{ij}=\tfrac14\omega^{is}\left( T^j_{nm;s}-2 R^j{}_{nms}\right)\extd x^m\wedge\extd x^n\in \Omega^2(M).\end{equation}

This solves Problem 1 in the introduction. Moreover, the classical connection $\nabla$ gets quantised as a map $\nabla_Q:\Omega^1\to \Omega^1\tens_1\Omega^1$ where $\tens_1$ means over the algebra
$A_\lambda$ with the above $\bullet_\lambda$ product. In noncommutative geometry quantum connections are handled like this as having values in an
extra copy of $\Omega^1$ which, classically, one would evaluate against a vector field to give the covariant derivative along that vector field.  The quantum
connection obeys two Leibniz rules
\[ \nabla_Q(a\bullet\xi)=a\bullet\nabla_Q\xi+\extd a\tens_1\xi,\quad \nabla_Q(\xi\bullet a)=(\nabla_Q\xi)\bullet a+\sigma_Q(\extd a\tens_1\xi)\]
where
\[ \sigma_Q:\Omega^1\tens_1\Omega^1\to \Omega^1\tens_1\Omega^1\]
is a bimodule map called the `generalised braiding' and is needed to make sense of the connection-like derivation property from the right.  If it exists we say that the quantum connection is a `bimodule connection' (and when it exists,
$\sigma_Q$ is unique so this is really a property of $\nabla_Q$ not extra data).

 In our case \cite{BegMa5} proves that there is such a bimodule connection
quantising our $\nabla$. In indices, it is
\[ \nabla_Q\extd x^i=-\left( \Gamma^i_{mn} +{\lambda\over 2}\omega^{sj}(\Gamma^i_{mk,s}\Gamma^k_{jn} -\Gamma^i_{kt}\Gamma^k_{sm}\Gamma^t_{jn} -\Gamma^i_{jk} R^k{}_{nms})\right)\extd x^m \tens_1 \extd x^n\]

In the quantisation process the torsion of $\nabla_Q$ gets modified to \cite{BegMa5}
\begin{eqnarray*}
T_{\nabla_Q}\xi={1\over 2}(\xi_i T^i_{nm}
+ \tfrac\lambda2\,(\partial_j  \righthalfcup  \nabla_i\xi)\,
\omega^{is}\,  T^j_{nm;s})\,\extd x^m\wedge_1\extd x^n\ .
\end{eqnarray*}

\subsection{Quantisation of the metric and Levi-Civita connection}

Here we suppose that $(M,\omega,\nabla)$ above has additional structure $(g,\widehat\nabla)$ where $g$ is a Riemannian (or pseudo-Riemannian) metric and $\widehat\nabla$ is the Levi-Civita connection. The condition in \cite{BegMa5} for the construction of a quantum metric comes out as the very reasonable requirement
\begin{equation}\label{qcen} \nabla g=0\end{equation}
that our Poisson connection is compatible with $g$.  It corresponds at the quantum level to the quantum metric being central and we assume the condition from now on. We build the quantum metric in two stages. First, the `functorial' choice is
\begin{equation}\label{gQ}
g_Q=q^{-1}_{\Omega^1,\Omega^1}(g)= g_{ij}\extd x^i\tens_1\extd x^j +{\lambda\over 2}\omega^{ij}g_{pm}\Gamma^p_{iq}\Gamma^q_{jn}\extd x^m\tens_1\extd x^n\end{equation}
where $q^{-1}$ will be explained in the next section but we have shown the result. One has $\nabla_Q g_Q=0$ so this is quantum metric compatible. We want our quantum metric to be `quantum symmetric' in the sense of killed by the quantum wedge product but we find 
\[ \wedge_1 g_Q=\lambda\CR;\quad \CR=H^{ij}g_{ij},\]
where explicitly
\[ \CR={1\over 2}\CR_{mn}\extd x^n\wedge\extd x^m,\quad \CR_{mn}={1\over 2}g_{ij}\omega^{is}(T^j_{nm;s}-R^j{}_{nms}- R^j{}_{mns})\ .\]
One can therefore either live with this or, which we do, define
\[ g_1=g_Q -{\lambda\over 4}g_{ij}\omega^{is}(T^j_{nm;s}-R^j{}_{nms}- R^j{}_{mns})\extd x^m\tens_1\extd x^n \]
as the quantum metric, which now has $\wedge_1(g_1)=0$. 

Finally, write the classical Levi-Civita connection as $\widehat\nabla=\nabla+S$ where
\[ S^i_{jk}=\tfrac12 g^{im}(T_{mjk}-T_{jkm}-T_{kjm})\]
is the contorsion tensor built from the torsion $T$ of $\nabla$. We do not want to give all the details but the main idea in \cite{BegMa5} is  to functorially quantise the contorsion to a quantum one $Q(S)$ and also to allow a further $O(\lambda)$ adjustment by a classical tensor $K$, i.e. we search for the quantum-Levi-Civita connection in the form
\[ \nabla_1=\nabla_{Q}+Q(S)+\lambda K.\]

\begin{theorem}\label{levithm}\cite{BegMa5} There is a unique quantum connection $\nabla_1$ which is quantum torsion free and for which the symmetric part in the last two factors of $\nabla_1g_1$ vanishes.  We have $\nabla_1g_1=0$ entirely if and only if
\[ \widehat\nabla \CR+ \omega^{ij}\,g_{rs}\,S^s_{jn}(R^r{}_{mki}+S^r_{km;i})\,\extd x^k\tens\extd x^m \wedge \extd x^n  =0\]
\end{theorem}
This is the condition (\ref{levicondition}) stated in tensor-calculus terms in the Introduction. The theorem says that there is always a unique `best-possible' quantum Levi-Civita connection but in general there could be an antisymmetric correction in the sense 
\[ (\id\tens\wedge)\nabla_1g_1=O(\lambda)\]
 as a new feature of quantum geometry. The condition for this correction to vanish is the one stated in the theorem but there are examples such as the axisymmetrically quantised Schwarzschild black hole in \cite{BegMa5} where we show that the correction cannot vanish for a large class of Poisson-connections $\nabla$. 

We will be interested only in the canonical special case  $S=T=0$ case of the above, where the $\nabla=\widehat\nabla$. In this case 
\[ \CR=-{1\over 2}g_{ij}\omega^{is} R^j{}_{nms}\extd x^m\wedge\extd x^n\]
and the `best possible' quantum Levi-Civita connection is just $\nabla_Q$ itself. The condition in the Theorem above is that this is covariantly constant. This holds for example for any K\"ahler-Einstein manifold and we will show the
results for $\mathbb{C P}^n$.

\subsection{Semiquantisation functor}

We do not want to explain all the mathematics here  but the main result in \cite{BegMa5} that then leads to the formulae above is as follows. We fix 
a manifold $M$ and a pair $(\omega,\nabla)$ of a Poisson tensor and Poisson-compatible connection. We show in \cite{BegMa5} that there
is a monoidal functor  
\[ Q:\{{\rm Vector\ bundles\ with\ connecton}\} \to \{{\rm Bimodules\ over}\ A_\lambda\ {\rm with\ connection}\}\]
This associates to a classical vector bundle $E$ and connection $\nabla_E$ on $M$ a bimodule $Q(E)$ over $A_\lambda$ equipped with a bimodule quantum connection $\nabla_{Q(E)}:Q(E)\to \Omega^1\tens_{A_\lambda}Q(E)$. Here $Q(E)$ is the sections of the classical bundle viewed with deformed left and right multiplication by $A_\lambda$:
\[ a\bullet_\lambda\xi=a\xi+{\lambda\over 2}\omega^{ij}a_{,i}\nabla_{Ej} \xi,\quad \xi\bullet_\lambda a=a\xi-{\lambda\over 2}\omega^{ij}a_{,i}\nabla_{Ej} \xi,\quad \forall a\in C^\infty(M),\ \xi\in Q(E).\]
Now both sides have a tensor product: the classical tensor product $\tens_0$ on the left means tensor product their sections, adding connections, on the right we have the usual tensor product $\tens_1$ of bimodules over and a certain tensor product of bimodule connections that makes use of their generalised braidings. The functor $Q$ being monoidal comes equipped with a functorial bimodule isomorphisms 
\[ q_{E,F}:(Q(E)\tens_1 Q(F))\to Q(E\tens_0 F),\quad q_{E,F}(\xi\tens_1\eta)=\xi\tens_0\eta+{\lambda\over 2}\omega^{ij}\nabla_{Ei}\xi\tens_0\nabla_{Fj}\eta\]
for all classical pairs $(E,\nabla_E), (F,\nabla_F)$, relating these tensor products. Formulae for the associated quantum bimodule connection are in \cite{BegMa5}. This functor is further extended to one where morphisms need not respect the connections, which is needed for the construction of $Q(S)$.

\section{Semiquantum Riemannian geometry of $\mathbb{C P}^n$}

Here we work out how the above general theory applies to $\mathbb{C P}^n$ as a K\"ahler-Einstein manifold. We first recall its classical geometry as real manifold so that we can directly use the formulae above. Full details will be in \cite{BegMa6} where we will  also discuss the noncommutative complex structure. 

\subsection{Geometry of $\mathbb{C P}^n$ as a Riemannian manifold}

We start with complex coordinates $w^i,\bar w^i$ for the sphere $S^{2n+1}\subset \R^{2n+2}$ with relations  
\begin{eqnarray*} \label{fubstu}
\sum w_i\,\bar w_i=1\ .\end{eqnarray*}
The quotient of $S^{2n+1}$ with its natural metric (inherited from twice the Euclidean metric) by the componentwise multiplicative action of $U(1)$ is the complex projective space $\mathbb{CP}^n$ with the Fubini-Study metric.  As usual we take coordinates $(z_1,\dots,z_n)$ for the open subset of $\mathbb{CP}^n$ where $w_0\neq 0$. This is done by setting
$(w_0,\dots,w_n)=(t,t\,z_1,\dots,t\,z_n)$ where
\begin{eqnarray*}
t^2=\frac1{1+|z_1|^2+\dots+|z_n|^2}\ .
\end{eqnarray*}
and finally we write $z^i=x^i+\mathrm{i}\, x^{i+n}$ for real coordinates $x^a$ where $1\le a\le 2n$. Outside this range it is convenient to use a `signed mod $2n$' rule where $x^b=-x^{b+2n}$ so that \begin{eqnarray*}
\tfrac{\partial x^c}{\partial x^a} &=& \kappa_{ac}:=\left\{\begin{array}{cc}+1 &\quad a-c\ \mathrm{an\ even\ multiple\ of}\ 2\,n \\-1 &\quad a-c\ \mathrm{an\ odd\ multiple\ of}\ 2\,n \\0 &\quad \mathrm{otherwise}\end{array}\right.\ ,\quad  \tfrac{\partial t}{\partial x^a} = -\,t^3\,x^a\ .
\end{eqnarray*}
In these coordinates the Fubini-Study metric, connection, curvature, symplectic and Poisson tensors are
\[ g^{ab}=\tfrac12\,t^{-2}\,(\kappa_{ab}+x^a\,x^b+x^{a+n}\,x^{b+n})\ ,\quad g_{ab}=2\,t^{2}\,\kappa_{ab}-2\,t^{4}(x^a\,x^b+x^{a+n}\,x^{b+n})\ . \]
\[ 
\Gamma^a_{bc} 
= 
 -\,  t^{2}\,  (x^c\,\kappa_{ab}  + x^b\,\kappa_{ac}  +
\kappa_{a+n,c}\,x^{b+n} + \kappa_{a+n,b}\,x^{c+n} )  \ . \]
\[ 
R^p{}_{cqb} 
      =   
\tfrac12\,g_{cb}  \, \kappa_{pq} -   \tfrac12\,g_{cq}   \, \kappa_{pb}   +   \tfrac12\,\omega_{bc}  \, \kappa_{p+n,q}     -    \tfrac12\,\omega_{qc}  \, \kappa_{p+n,b}     +      \kappa_{p+n,c}\, \omega_{bq}  \ .
\]
\[ \omega_{ab}=g_{a,b+n}=2\,t^2 \kappa_{a,b+n}-2\,t^4(x^a x^{b+n}-x^{a+n}x^b),\quad 
\omega^{ab}=\tfrac12\,t^{-2}\,(\kappa_{a,b+n}+x^a\,x^{b+n}-x^{a+n}\,x^b)\ . \]
The curvature map in our conventions is 
\[  R_\nabla= {1\over 2}\extd x^a\wedge \extd x^b\tens[\nabla_a,\nabla_b],\quad [\nabla_a,\nabla_b]\extd x^c=-R^c{}_{dab}\extd x^d\ \]
and the symplectic 2-form (in our conventions) is 
\[ \varpi=\omega_{ab}\extd x^b\wedge\extd x^a\ .\]

Using the formula (\ref{H})  gives
\begin{align*}
H^{ab}&=\tfrac14\,\extd x^{a+n}\wedge \extd x^b+\tfrac14\, \extd x^{b+n}\wedge\extd x^a-  \tfrac14\,  g^{ab}\varpi  \ .
\end{align*}
\[ \CR=-\tfrac12\,(n+1)\varpi\ \]
so that the Ricci 2-form is a multiple of the symplectic 2-form, just as the usual Ricci curvature is a multiple of $g$. 

If we want some of these in complex coordinates then we have 
\[ g=g_{i \, \bar j}\,(
\extd z^i\tens \extd \bar z^j + \extd \bar z^j\tens \extd z^i)\ ,\quad g_{i\,\bar j}=(t^2\,\delta_{ij}-t^4\, \bar z^i\,z^j)\ .\]

It will be convenient to define a 1-form 
\[ \tau={\sum_{i=1}^n \bar z^i\extd z^i\over 1+ |\vec z|^2}\]
in our patch. Here $\tau+\bar\tau =\extd\ln(1+|\vec z|^2)=-t^{-2}\extd t^2$. It will also be convenient for notation to set $z^i_+=z^i$ and 
$z^i_-=\bar z^i$, and correspondingly $\tau_+=\tau$ and $\tau_-=\bar\tau $. 

We have followed standard conventions in defining the canonical 1-form for a K\"ahler manifold by $\varpi=(J\wedge \id)g$ which gives
\begin{align*}
\varpi=2 \mathrm{i} g_{i\bar j}\extd z^i\wedge\extd \bar z^j\  = -2\mathrm{i} \extd\tau\ .
\end{align*}

\medskip
We may also recall that $\mathbb{C P}^n$ is an example of a K\"ahler manifold. Hence its structure is given via a 
K\"ahler potential which, in the case of $\mathbb{C P}^n$ is
\[ K_j=\ln\Big(\sum_{a=0}^{n}|{w^a\over w^j}|^2\Big)\]
in the coordinate chart $U_j=\{(w^0,\cdots,w^{n})\ |\ w^j\ne 0\}$. 
On $U_0$ with our complex coordinates $z_1,\dots,z_n$, $K_0=\ln(1+|\vec z|^2)$. Then we calculate 
 $\tau=\partial K_0$ and
 \begin{align*}
\frac{\partial^2 K_0}{\partial z_i\,\partial \bar z_j}\,=\, g_{i\bar j}\ , \quad \varpi=2\,\mathrm{i}\,\partial\pdol K_0\ .
\end{align*}

Finally, we introduce \[ 
 \gamma_+= \gamma={t^2}\extd \bar z^i\tens\extd z^i-\bar\tau \tens\tau,\quad  \gamma_-=\bar\gamma={t^2}\extd z^i\tens\extd \bar z^i-\tau\tens\bar\tau\]
with summation understood. Then
\[ g=\gamma+\bar\gamma,\quad  \varpi=\mathrm{i}\wedge(\bar\gamma-\gamma)=-2 \mathrm{i}\wedge(\gamma).\]

\begin{proposition}\label{levi} For $\mathbb{C P}^n$, the Levi-Civita connection and its curvature associated to the standard complex structure and the Fubini-Study metric are
\[  \nabla \extd z^i_\pm=\tau_\pm\tens\extd z^i_\pm+\extd z^i_\pm\tens\tau_\pm\ ,\quad  R_\nabla(\extd z_\pm)=\pm { \mathrm{i}\over 2}\varpi\tens\extd z^i_\pm -\extd z^i_\pm \wedge \gamma_\pm\]
\end{proposition}

The curvature map here was computed using algebraic methods but one can check that it agrees with the tensor calculus computation\cite{BegMa6}. We have not seen such a simple description of the geometry of $\mathbb{C P}^n$ elsewhere but presumably this is known. There are similar formulae in other coordinate patches.

\subsection{Semiquantum Riemannian geometry of  $\mathbb{C P}^n$}
According to our constructions we have the  following quantised product, for functions $e,f$:
\begin{align}
e\bullet f= e\,f + \frac{\mathrm{i}\,\lambda}{2}\, t^{-2}\,(\delta_{ij}
+ z^i\, \bar z^j)(\pd_i e\  \pdol_j f    - \pdol_j  e \ \pd_i f)\ .
\end{align}
In fact this formula also gives the quantised product of a function and a form, if we replace one of $e$ or $f$ by a form and use the (complex) Levi-Civita
connection for $\pd_i$ and $\pdol_j$. The undeformed  cases are
\begin{align*}
z^i\bullet z^j=z^i\, z^j\ ,\quad z^i\bullet \extd z^j=z^i\, \extd z^j\ ,\quad \extd z^i\bullet z^j=\extd z^i\, z^j\ ,
\end{align*}
and the same formulae hold if we bar all the $z$s. As the exterior derivative is undeformed, applying $\extd$ to these results gives
$\extd z^i\wedge_1 \extd z^j=\extd z^i\wedge \extd z^j$ and $\extd \bar z^i\wedge_1 \extd \bar z^j=\extd \bar z^i\wedge \extd \bar z^j$. 
However when we mix $z$s and $\bar z$s in the same product, we get non commutative behaviour:
\begin{align*}
z^i\bullet \bar z^j=z_i\, \bar z^j+ \frac{\mathrm{i}\,\lambda}{2}\, t^{-2}\,(\delta_{ij}
+ z^i\, \bar z^j)\ ,\quad \bar z^j\bullet  z^i=z_i\, \bar z^j- \frac{\mathrm{i}\,\lambda}{2}\, t^{-2}\,(\delta_{ij}
+ z^i\, \bar z^j)\ ,
\end{align*}
so we can write a commutation relation
\begin{align*}
[z^i,z^j]_\bullet=0=[\bar z^i,\bar z^j]_\bullet,\quad 
[z^i, \bar z^j]_\bullet= {\mathrm{i}\,\lambda}\, t^{-2}\,(\delta_{ij}
+ z^i\, \bar z^j)\ .
\end{align*}
If we mix functions and forms, we get,
\begin{align*} \label{bulletone}
z^i\bullet \extd \bar z^j=&\ z^i\, \extd \bar z^j + \tfrac{\mathrm{i}\,\lambda}{2}\, t^{-2}\,(\delta_{ij}
+ z^i\, \bar z^j)\, \bar\tau  
 + \tfrac{\mathrm{i}\,\lambda}{2}\,t^{-2}\, z^i\, \extd \bar z^j  \ ,\cr
 \extd \bar z^j \bullet z^i=&\ z^i\, \extd \bar z^j - \tfrac{\mathrm{i}\,\lambda}{2}\, t^{-2}\,(\delta_{ij}
+ z^i\, \bar z^j)\, \bar\tau  
 - \tfrac{\mathrm{i}\,\lambda}{2}\,t^{-2}\, z^i\, \extd \bar z^j \ .
\end{align*}
which gives the commutation relations
\[[z^i,\extd z^j]_\bullet=0=[\bar z^i,\extd\bar z^j]\]
\[  [z^i, \extd \bar z^j]_\bullet=\mathrm{i}\,\lambda\,  t^{-2}\, \left((\delta_{ij}
+ z^i\, \bar z^j)\, \bar\tau  
 +  z^i\, \extd \bar z^j\right)\ , \quad [\bar z^i, \extd z^j]_\bullet=-\mathrm{i}\,\lambda\,  t^{-2}\, \left((\delta_{ij}
+ \bar z^i\, z^j)\, \tau  
 +  \bar z^i\, \extd  z^j\right)\ .\]

For the wedge product using Section~2.1 we have 
\begin{align*} \label{bullettwo}
\extd z^i\wedge_1 \extd \bar z^j=&\ \extd z^i\wedge \extd \bar z^j +  \tfrac{\mathrm{i}\,\lambda}{2}\, t^{-2}\, \big((\delta_{ij}
+ z^i\, \bar z^j)\, t^2\,\extd z^k\wedge\extd \bar z^k  + \tau \wedge z^i\, \extd \bar z^j +
\bar z^j \extd z^i\, \wedge \bar\tau   + \extd z^i\wedge \extd \bar z^j  \big)\ ,\cr
\{\extd z^i, \extd \bar z^j\}_{\wedge_1}=&\ \mathrm{i}\,\lambda\, t^{-2}\, \big((\delta_{ij}
+ z^i\, \bar z^j)\, t^2\,\extd z^k\wedge\extd \bar z^k  + \tau \wedge z^i\, \extd \bar z^j +
 \bar z^j\extd z^i\, \wedge \bar\tau   + \extd z^i\wedge \extd \bar z^j  \big)\ .
\end{align*}

Let us note that to order $\lambda$ we are free to write these as some kind of $q$-commutator where
\[ q=e^{\mathrm{i} \lambda t^{-2}}\]
\[ q \bar z^i z^j-z^j \bar z^i={\lambda t^{-2}\over \mathrm{i}}\delta_{ij}\]
\[ q \bar z^i \extd z^j-(\extd z^j)\bar z^i={\lambda t^{-2}\over \mathrm{i}}(\delta_{ij}+\bar z^iz^j)\tau\]
\[ q^{-1} z^i \extd \bar z^j-(\extd \bar z^j)z^i=-{\lambda t^{-2}\over \mathrm{i}}(\delta_{ij}+z^i\bar z^j)\bar\tau \]
\[ q^{-1}\extd z^i\wedge_1\extd \bar z^j+\extd\bar z^j\wedge_1\extd z^i=\ \mathrm{i}\,{\lambda}\, t^{-2}\, \big((\delta_{ij}
+ z^i\, \bar z^j)\, t^2\,\extd z^k\wedge\extd \bar z^k  + \tau \wedge z^i\, \extd \bar z^j +\bar z^j\extd z^i\, \wedge \bar\tau  \big)\ .\]

We can also compute these relations on our constrained homogeneous coordinates $w^i$ where $w^0=t$ is real and positive and $w \bar w=1$ as $n+1$-vectors (i.e. $t^2=1-|\vec w|^2$ from this point of view of $w^i=t z^i$ as the complex coordinates).  

\begin{proposition} In the `upstairs' coordinates restricted to $i,j>0$ we have
\[ [w^i,w^j]_\bullet=0,\quad [w^i,\bar w^{j}]_\bullet= \mathrm{i}\lambda\delta_{ij}\]
\[ [w^i,\extd \bar w^{j}]_\bullet=\mathrm{i}{\lambda\over 2} \left((2\delta_{ij}+w^i\bar w^{j}{(t^{-2}-2)}){(\bar\tau -\tau)\over 2}+ w^i\extd \bar w^{j}-\bar w^{j}\extd w^i\right)\]
\[ [w^i,\extd w^j]_\bullet={\lambda\over 2 \mathrm{i}}\left(w^iw^j(2\bar\tau +t^{-2}{(\bar\tau -\tau)\over 2})+w^i\extd w^j+w^j\extd w^i\right)\]
\end{proposition}

We see that our algebra relations agree with the recent proposal\cite{Pen} of Roger Penrose for `quantum twistor space' but in our case with the Euclidean signature. On the other hand, these restricted homogeneous coordinates are less well-adapted to the
holomorphic nature of the calculus even if they put the algebra commutation relations in canonical form. Whichever coordinates are used, it should be remembered that while the coordinate algebra can be constructed to all orders in $\lambda$ associatively (for example using geometric quantisation via $\tau$ as a connection with curvature yielding $\varpi$) this is not the case for the differential calculus which, since $\nabla$ on $\mathbb{C P}^n$ has curvature, will be nonassociative at order $\lambda^2$.

Finally, a long computation in \cite{BegMa6} but using the general results in Section~2.1 gives us the quantum metric and quantum-Levi-Civita connection. Because we have taken
$\nabla$ to be the Levi-Civita connection and because $\CR$ is a multiple of the metric it is covariantly constant,  Theorem~\ref{levithm}
in Section~2.2 applies but in the simplified form where we the quantum Levi-Civita connection is just $\nabla_Q$ itself. 

\begin{proposition}
\begin{align*}
g_1 = g_{i \, \bar j}\,\extd z^i\tens_1 \extd \bar z^j + g_{i \, \bar j}\,\extd \bar z^j\tens_1 \extd z^i + \tfrac{\lambda}2\,(n+1)\, \widetilde\varpi ,
\end{align*}
where $\widetilde\varpi=\mathrm{i}(\bar\gamma-\gamma)$ is taken with $\tens_1$, and
\begin{align*}
 \nabla_Q \extd z^i_\pm =& \  (1\pm\mathrm{i}\,\lambda)\,( \tau_\pm\tens_1\extd z^i_\pm+\extd z^i_\pm\tens_1\tau_\pm) \ .
\end{align*}
\end{proposition}

We note that the quantum $\nabla_Q$ has a strikingly similar form to the classical Levi-Civita connection in Theorem~\ref{levi}. The general theory means that there is an associated $\sigma_Q$ making it a bimodule connection.

\section{Quantum geometry in quantum mechanics}

The currently envisaged application of the above is with $\lambda$  related to the Planck scale, i.e. applying the theory to quantum spacetime.
 However, here we want to consider the question of ordinary
quantum mechanics where $\lambda$ should be related to $\hbar$ and other physical scales in ordinary quantum systems.

Specifically, the questions we pose, if $(M,\omega)$ is a quantum mechanical phase space, are:

\begin{quote} {\bf Question 1} What is the physical content of quantum differential forms on $A_\lambda$? \end{quote}

\begin{quote} {\bf Question 2} What is the physical content of  quantum metrics and connections on the quantisation on $A_\lambda$? \end{quote}

Let's consider the first question. In classical mechanics the phase space is not merely a topological space, it is a manifold and this differential structure
is used in formulating the Hamilton-Jacobi equations of motion
\[ \dot a= \{a,H\}=-\hat H(a),\quad \forall a\in C^\infty(M).\]
In other words, time evolution is by the vector field $-\hat H$. When we quantise we might then expect the quantum evolution to be given to lowest order by a quantum
version of the Hamiltonian-Jacobi equations using the quantum differential forms. Of course, the usual proposal is to replace
Poisson bracket by commutator:
\[ \dot a={\mathrm{i}\over\hbar}[H,a]=\{a,H\}+O(\hbar)\]
while we might in quantum geometry be inclined to  something like
\[ \dot a=\omega_1(\extd a\tens_1\extd H)\]
for a natural quantum Poisson tensor  $\omega_1$ and a quantum pairing, to be constructed. This comparison would give us partial information about the next (2nd) order terms in  product of $A_\lambda$. We can also compare with this order in the Fedosov quantisation $A_\lambda$ which is determined through a flat symplectic connection that we could also use as $\nabla$. This issue remains further to be investigated. 

We also should consider the question of {\em how do differential forms evolve in quantum mechanics?} If we take the view that they do so by commutator with
$H$ then 
\[ \dot{(\extd a)}={\mathrm{i}\over\hbar}[H,\extd a]=-\nabla_{\hat H}(\extd a)+O(\hbar)\]
 which is very reasonable for it says that we use the same vector field $-\hat H$ but now with the covariant derivative given by our Poisson connection. Note that this is not the Lie derivative so
time evolution is not a diffeomorphism of the classical system. Rather, our induced classical picture is that as a particle moves along a trajectory with tangent vector given by $-\hat H$, any differentials are parallel transported also using the connection $\nabla$. This should also apply to the evolution of points in other bundles over phase space that are equipped with connections as in Section~2.3. On the other hand, by the Poisson-compatibility condition (\ref{poicomp}), the above classical evolution of 1-forms is equivalent to
\[\dot{(\extd a)} =\extd\dot a-\nabla_{\hat a}(\extd H)\]
which reminds us that time evolution does {\em not} commute with $\extd$ unless the 2nd term vanishes. As for functions, we can ask about the $O(\hbar)$ corrections to this evolution and compare with the quantum geometry via $\nabla_Q$.

To put some of these ideas in concrete terms, lets look at the simplest case $M=\R^{2n}$ and canonical coordinates $\{q^i,p^i\}$, Euclidean metric and trivial $\nabla$ so that
\[ [q^i,p^j]=\mathrm{i}\hbar\delta_{ij},\quad g=\extd q^i\tens\extd q^i+\extd p^i\tens\extd p^i,\quad \nabla\extd q^i=\nabla\extd p^i=0.\]
Here $\dot{(\extd q^i)}=\dot{(\extd p^i)}=0$ which means that our choice of coordinate basis for $\Omega^1$ is not affected by time evolution. On the other hand, both $\dot q^i,\dot p^i$ are
not normally constant on $M$ as they are given by the Hamiltonian-Jacobi equations. For example
\[ H={p^2\over 2m}+V(q)\quad \Rightarrow\quad \dot{(\extd a)}-\extd\dot a=-{1\over m}{\del a\over\del q^i}\extd p^i+ {\del^2 V\over\del q^i\del q^j}{\del a\over\del p^j}\extd q^i\ .\]
Thus our proposal seems reasonable for the evolution of differential forms, but it is still an assumption that should be put to experimental test. Or rather, such an evolution may be natural in classical  mechanics when we have internal geometric structure at each point of phase space. The quantum calculus, meanwhile, has the same form on the generators as classically and is associative as $\nabla$ is trivial and flat. 

Question 2 about the quantum Riemannian geometry of phase space entails a prequestion  about the classical Riemannian geometry of phase space. In physics, one place where this enters is in the description of Berry phase. For example on $\mathbb{C P}^n$ seen as a state space of a quantum system its Riemannian geometry enters into this and into expressions for higher-power uncertainty relations\cite{Brody}. In the K\"ahler case such as this, the metric is canonical given the symplectic and complex structures but in other cases the prequestion is what should be the physical significance of the metric on phase space? Also note that we have different choices for $\nabla$ and if we take the Levi-Civita connection we will tend to have nonassociativity of the calculus in the presence of Riemannian curvature as in our $\mathbb{C P}^n$ example. However, when the manifold is parallelizable, one can also take the Weitzenb\"ock $\nabla$ as in teleparallel gravity\cite{tele}, which is flat but has torsion. Then we will need the general case of Theorem~\ref{levithm}.

There has also been much interest recently in an interpretation\cite{Ma:mean} of noncommutative spacetime as {\em curved} momentum space or `cogravity' in some sense. In the same
way by {\em quantum Born reciprocity}, a curved position spacetime should correspond locally to noncommutative position space. Thus at the Poisson level non-zero $\omega$ in the $q$ sector of phase space should correspond to cogravity in the spatial momentum while non-zero $\omega$ in the $p$ sector should be a signal of gravity or at least of curvature in space. Hence, when there is quantum noncommutativity of space or spacetime {\em and} gravity, there should be both $\omega$ and phase space curvature generically, which is how quantum groups first arose out of Planck scale physics in \cite{Ma:pla}. The problem in this context is what precise equations should govern the interaction of these structures and what (\ref{levicondition}) provides is a first instance of such an equation, coming from the assumption of the existence of a quantum geometry of which the classical manifold is a classical limit.  In this context a natural special case of (\ref{levicondition}) would be where $\omega$ is zero among half the variables and the curvature zero among the other half (the gravity or cogravity special cases). Another idea is {\em quantum Born reciprocity} that one should be able to interchange position and momentum as a key idea self-duality principle for quantum gravity\cite{Ma:pla} and this could now be imposed at the semiclassical level as a further condition on (\ref{levicondition}). We do not assume that $\nabla$ is the Levi-Civita connection but instead we need new principles and equations to help determine it.

\end{document}